\documentclass[]{interact}

\usepackage{epstopdf}

\usepackage{natbib}
\bibpunct[, ]{[}{]}{,}{n}{,}{,}
\renewcommand\bibfont{\fontsize{10}{12}\selectfont}

\usepackage{hyperref}
\numberwithin{equation}{section}

\usepackage{wrapfig}

\RequirePackage{amsthm,amsmath,amsfonts,amssymb}

\usepackage{wrapfig,lipsum,booktabs}
\usepackage{caption}
\usepackage{graphicx}
\usepackage{subcaption}
\usepackage{placeins}
\usepackage{caption}

\theoremstyle{plain}

\theoremstyle{definition}

\theoremstyle{remark}

\begin{document}
	
	\articletype{ARTICLE TEMPLATE}
	
	\title{A Note On Simultaneous Estimation  of Order Restricted Location Parameters of a General Bivariate Symmetric Distribution Under a General Loss Function}
	
	\author{
		\name{Naresh Garg\textsuperscript{a}\thanks{CONTACT Naresh Garg. Email: nareshng@iitk.ac.in} and Neeraj Misra\textsuperscript{a}}
		\affil{\textsuperscript{a}Department of Mathematics and Statistics, Indian Institute of Technology Kanpur, Uttar Pradesh, India}
	}
	
	\maketitle
	
	\begin{abstract}
		The problem of simultaneous estimation of order restricted location parameters $\theta_1$ and $\theta_2$ ($-\infty<\theta_1\leq \theta_2<\infty$) of a bivariate location symmetric distribution, under a general loss function, is being considered. In the literature, many authors have studied this problem for specific probability models and specific loss functions. In this paper, we unify these results by considering a general bivariate symmetric model and a quite general loss function. We use the Stein and the Kubokawa (or IERD) techniques to derive improved estimators over any location equivariant estimator under a general loss function. We see that the improved Stein type estimator is robust with respect to the choice of a bivariate symmetric distribution and the loss function, as it only requires the loss function to satisfy some generic conditions. A simulation study is carried out to validate the findings of the paper. A real-life data analysis is also provided.
	\end{abstract}
	

	\section{Introduction}\label{sec1}
	
	The problem of estimating order restricted location parameters $\theta_1$ and $\theta_2$ ($-\infty<\theta_1\leq \theta_2<\infty$) of two distributions is of interest in many real-life situations. For example, in an engine efficiency measurement experiment where estimating the average efficiency of an internal combustion engine (IC engine) and an external combustion engine (EC engine) is of interest, it can be assumed that the average efficiency of an IC engine is higher than the average efficiency of an EC engine. For an account of such applications and relevant literature, one may refer to Barlow et al. (\citeyear{MR0326887}), Robertson et al. (\citeyear{MR961262}), and van Eeden (\citeyear{MR2265239}).
	\vspace*{2mm}

	Early studies in this area were focused on studying isotonic regression and/or restricted maximum likelihood estimators of order restricted parameters. Afterwards, the problem was studied using decision theoretic approach with focus on obtaining estimators improving over the unrestricted best location equivariant estimators (BLEE) and/or unrestricted maximum likelihood estimators (MLE). Several of these studies are centered around specific distributions and specific loss functions, barring a few studies that are carried out for general probability models and general loss function. In this paper, we will obtain some unified results for simultaneous estimation of order restricted location parameters of a general bivariate symmetric distribution under a general loss function. \vspace{2mm}
	
	It is worth mentioning that Stein (\citeyear{MR171344}) proposed a technique to improve the best affine equivariant estimator of the variance of a normal distribution. The Stein technique provides shrinkage type non-smooth dominating estimators. This technique is generalized by Brewster and Zidek (\citeyear{MR381098}) who also proposed another technique to improve the best equivariant estimators. The Brewster and Zidek (\citeyear{MR381098}) technique produces smooth dominating estimators (generally the generalized Bayes estimator with respect to a non-informative prior). Kubokawa (\citeyear{MR1272084}) unified the two techniques of the Stein (\citeyear{MR171344}) and the Brewster and Zidek (\citeyear{MR381098}) by using a representation of difference of risk functions in terms of a definite integral. He named this unified technique the integral expression risk difference (IERD) method. \vspace{2mm}
	
	Using the techniques of Stein (\citeyear{MR171344}), Brewster and Zidek (\citeyear{MR381098}) and Kubokawa (\citeyear{MR1272084}), several authors have considered the problem of improving equivariant estimators for specific probability models (mostly, having independent marginals) and specific loss functions. Kumar and Sharma (\citeyear{MR981031}) dealt with simultaneous estimation of ordered means of two normal distributions, having a known common variance, under the sum of squared error loss functions and obtained a sufficient condition that ensures the inadmissibility of any location equivariant estimator. Patra and Kumar (\citeyear{patra2017estimating}) further extended the results of Kumar and Sharma (\citeyear{MR981031}) to a bivariate normal distribution having known variances and a known correlation coefficient. Tsukuma and Kubokawa (\citeyear{tsukuma2008stein}) considered simultaneous estimation of $p \;(\geq 2)$ means of a $p$ dimensional multivariate normal distribution with the covariance matrix as the identity matrix, when it was known apriori that the means were restricted to a polyhedral convex cone. Under the sum of the squared errors loss function, they have obtained the generalized Bayes estimator against the uniform prior distribution over the polyhedral convex cone and shown that this estimator is minimax. Further, for $p=2$ and when means are order restricted, they have shown that the generalized Bayes estimator is admissible.  For a general framework, Hamura and Kubokawa (\citeyear{hamura2022robustness}) considered component-wise estimation of order restricted location parameters of two independent log-concave or log-convex probability models. They have considered the Stein type truncated estimator and shown that this estimator dominates the usual unrestricted estimator, under the squared error loss function.

	\vspace*{3mm}

	In this paper, we aim to unify various results in the literature by considering a general bivariate symmetric location model and a general loss function. We consider simultaneous estimation of order restricted location parameters $\theta_1$ and $\theta_2$ ($\theta_1\leq \theta_2$), under a general loss function. We use the Stein (1964) technique to derive a sufficient condition for the inadmissibility of a location equivariant estimator and obtain an improved estimator. In this case, the improved Stein (1964) type estimator is robust, as the form of the improved estimator does not depend on the choice of the bivariate symmetric distribution and the loss function, except some generic conditions. Further, we use the Kubokawa (1994) technique to obtain a class of improved estimators over the best location equivariant estimators (BLEE). To illustrate the usefulness of our results, we consider a bivariate normal distribution with unknown order restricted means, the known common variance, and the known correlation coefficient. We obtain the Stein (1964) type improved estimator over the unrestricted BLEE/MLE under a general loss function. We see that this improved estimator is the restricted MLE. We also obtain the Brewster-Zidek (1974) type improved estimators over the unrestricted BLEE/MLE under the squared error loss and the absolute error loss. These improved estimators are also generalized Bayes estimators under the squared error loss and the absolute error loss, respectively.
	\vspace*{3mm}

	The rest of the paper is organised as follows: In Section \ref{sec2}, we consider the estimation of location parameters of a general bivariate symmetric model. In Section \ref{sec2.1}, we use the Stein (1964) technique to show the inadmissibility of location equivariant estimators satisfying a sufficient condition, and in Section \ref{sec2.2}, we use the Kubokawa (1964) technique to obtain a class of estimators improving upon the BLEE, under a general loss function. In Section \ref{sec3}, we demonstrate applications of our general result to a bivariate normal distribution and report a simulation study to validate our findings. In Section \ref{sec4}, we present a real-life application of the findings of this paper. In Section \ref{sec5}, we provide concluding remarks for this paper.

	\section{\textbf{Improved estimators for order restricted location parameters}}
	\label{sec2}
	
	Let $\bold{X}=(X_1,X_2)$ be a random vector having the Lebesgue probability density function (p.d.f.) 
	\begin{equation}\label{eq:2.1}
		f_{\boldsymbol{\theta}}(x_1,x_2)=f(x_1-\theta_1,x_2-\theta_2),\; (x_1,x_2)\in \Re^2, \;\boldsymbol{\theta}=(\theta_1,\theta_2)\in \Theta_0,
	\end{equation}
	where $f(\cdot,\cdot) $ is a specified bivariate Lebesgue p.d.f. on $\Re^2=(-\infty,\infty)\times(-\infty,\infty)$, $\boldsymbol{\theta}=(\theta_1,\theta_2) \;(\in \Theta_0)$ is an unknown parameter and $\Theta_0$ is the parameter space. Generally, $\bold{X}=(X_1,X_2)$ would be a minimal-sufficient statistic for $\boldsymbol{\theta}\in \Theta_0$, based on a bivariate random sample or two independent random samples, as the case may be. Throughout, we make the following assumption about the probability model \eqref{eq:2.1}:
	\\~\\ \textbf{Assumption D1:} The parameter space of interest is the restricted space $\Theta_0=\{(x,y):-\infty<x\leq y<\infty\}. \text{ Moreover,}\; f(z_1,z_2)=f(z_2,z_1)=f(-z_1,-z_2),\; \forall\; (z_1,z_2)\in \Re^2$.

	\vspace*{2mm}
	
	\noindent Consider simultaneous estimation of order restricted location parameters $\theta_1$ and $\theta_2$ ($(\theta_1,\theta_2)\in \Theta_0$) under the loss function
	\begin{equation}\label{eq:2.2}
		L(\boldsymbol{\theta},\bold{a})= W(a_1-\theta_1) + W(a_2-\theta_2),\;\boldsymbol{\theta}=(\theta_1,\theta_2)\in\Theta_0,\;\bold{a}=(a_1,a_2)\in\mathcal{A}=\Re^2,
	\end{equation}
	where $W:\Re\rightarrow [0,\infty)$ satisfies the following assumption:
	\\~\\\textbf{Assumption D2:}  $W(0)=0$, $W(t)=W(-t),\; t\in \Re$, $W(t)$ is strictly decreasing in $t\in (-\infty,0)$ and strictly increasing in $t\in (0,\infty)$. Also $W^{'}(t)$ is increasing, almost everywhere.
	\vspace*{2mm}

	The above estimation problem is invariant under the group $\mathcal{G}=\{g_c:\,c\in\Re\}$ of transformations, where $g_c(x_1,x_2)=(x_1+c,x_2+c),\; (x_1,x_2)\in\Re^2,\;c\in\Re.$ The induced group of transformations on the parameter space $\Theta_0$ and the action space $\mathcal{A}$ are $\overline{\mathcal{G}}=\{\overline{g}_{c}:c\in\Re \}$ and $\tilde{\mathcal{G}}=\{\tilde{g}_{c}:c\in\Re \}$, respectively, where $\overline{g}_{c}(\theta_1,\theta_2)=(\theta_1+c,\theta_2+c),\;(\theta_1,\theta_2)\in\Theta_0,$ $\tilde{g}_{c}(a_1,a_2)=(a_1+c,a_2+c)$, $\;(a_1,a_2)\in\mathcal{A}=\Re^2,\;  c\in\Re$. Any location equivariant estimator of $\boldsymbol{\theta}$ is of the form
	\begin{equation}\label{eq:2.3}
		\boldsymbol{\delta}_{\psi}(\bold{X})=(X_1-\psi_1(D),X_2-\psi_2(D)),
	\end{equation}
	for some functions $\psi_i:\,\Re\rightarrow \Re,\;i=1,2$, where $D=X_2-X_1$.
	\vspace*{2mm}

	Let $Z_i=X_i-\theta_i,\;i=1,2$. Since $f(z_1,z_2)=f(z_2,z_1),\; \forall \; (z_1,z_2)\in \Re^2$ and $f(z_1,z_2)=f(-z_1,-z_2),\; \forall \; (z_1,z_2)\in \Re^2$, we have $(Z_1,Z_2)\overset{\mathrm{d}}{=}(Z_2,Z_1)\overset{\mathrm{d}}{=}(-Z_2,-Z_1),$ where $\overset{\mathrm{d}}{=}$ stands for equality in the distribution. Evidently, the problem of simultaneously estimating the order restricted location parameters $\theta_1$ and $\theta_2$ ($\boldsymbol{\theta}\in \Theta_0$), under the loss function (2.2), is also invariant under the group of transformations $\mathcal{H}=\{h_1,h_2\},$ where $h_1(x_1,x_2)=(x_1,x_2),\; h_2(x_1,x_2)=(-x_2,-x_1),\; (x_1,x_2)\in \Re^2.$ The induced group of transformations on the parameter space $\Theta_0$ and the action space $\mathcal{A}$ are $\overline{H}$ and $\tilde{H}$, respectively, where $\overline{H}=\{\overline{h}_1,\overline{h}_2\},\;  \tilde{H}=\{\tilde{h}_1,\tilde{h}_2\},\; \overline{h}_1(\theta_1,\theta_2)=(\theta_1,\theta_2),\; \overline{h}_2(\theta_1,\theta_2)=(-\theta_2,-\theta_1),\; (\theta_1,\theta_2)\in \Theta_0,\; \tilde{h}_1(a_1,a_2)=(a_1,a_2)$ and $\tilde{h}_2(a_1,a_2)=(-a_2,-a_1),\; (a_1,a_2)\in \mathcal{A}$. An estimator $\bold{\delta}_{\psi}(\bold{X})=(X_1-\psi_1(D),X_2-\psi_2(D))$ is invariant under $\mathcal{H}$ if, and only if, 
	$$(X_1-\psi_1(D),X_2-\psi_2(D))=(-(-X_1-\psi_2(D)),-(-X_2-\psi_1(D)))$$
	i.e., $\psi_2(D)=-\psi_1(D)$.
	\vspace*{2mm}
	
	\noindent Thus, the form of any estimator that is equivariant under $\mathcal{G}$ as well as under $\mathcal{H}$ is
	\begin{equation}\label{eq:2.4}
		\boldsymbol{\delta}_{\psi}(\bold{X})=(X_1-\psi(D),X_2+\psi(D)),
	\end{equation}
	for some function $\psi:\Re\rightarrow \Re$.

	\subsection{Improvements Over an Location Equivariant Estimator $(X_1-\psi(D),X_2+\psi(D))$}\label{sec2.1}

	\vspace*{2mm}
	
	\noindent
	
	\vspace*{4mm}
	
	In this section, we use the Stein (1964) technique to obtain improved estimators over any arbitrary location equivariant estimator $\boldsymbol{\delta}_{\psi}(\bold{X})=(X_1-\psi(D),X_2+\psi(D))$. Using $f(z_1,z_2)=f(-z_1,z_2),\; (z_1,z_2)\in \Re^2$ (i.e. $(Z_1,Z_2)\overset{\mathrm{d}}{=} (-Z_2,-Z_1)$), and $W(t)=W(-t),\; t\in \Re$, for $\boldsymbol{\theta}=(\theta_1,\theta_2)\in \Theta_0$, the risk function of the estimator $\delta_{\psi}(\bold{X})$, defined by \eqref{eq:2.4}, is obtained as
	\begin{align*}
		R(\boldsymbol{\theta},\delta_{\psi})&=E_{\boldsymbol{\theta}}[W(X_1-\psi(D)-\theta_1)+W(X_2+\psi(D)-\theta_2)]\\
		&=2\int_{-\infty}^{\infty} \left[\int_{-\infty}^{\infty} W(s-\psi(t)) f(s,s+t-\lambda)ds \right] dt\\
		&=2 \int_{-\infty}^{\infty} r_{\lambda}(\psi(t),t)\,dt,
	\end{align*}
	where $\lambda=\theta_2-\theta_1\; (\geq \,0)$ and
	$$r_{\lambda}(c,t)=\int_{-\infty}^{\infty} W(s-c) f(s,s+t-\lambda) ds,\;\;c\in \Re,\; t\in \Re.$$
	\vspace*{2mm}
	
	\noindent The following lemma will be useful in proving the main result of the paper.
	\vspace*{2mm}

	\noindent \textbf{Lemma 2.1.1.} Suppose that the assumptions (D1) and (D2) hold. For any $t\in \Re$, $\lambda \geq 0$, and $c\in \Re$,
	$$ r_{\lambda}(c,t)\geq r_{\lambda}\left(\frac{\lambda-t}{2},t\right).$$
	\begin{proof}
		Let $t\in \Re$, $\lambda \geq 0$, and $c\in \Re$ be fixed. Then
		$$r_{\lambda}(c,t)=\int_{-\infty}^{\infty} W(s-d) f\left(s+\frac{\lambda-t}{2},s-\frac{\lambda-t}{2}\right) ds=r^{*}_{\lambda}(d,t), \text{ say},$$
		where $d=c-\frac{\lambda-t}{2}$. We have
		\begin{align*}
			r^{*}_{\lambda}(d,t)&=\int_{-\infty}^{\infty} W(s-d) f\left(s+\frac{\lambda-t}{2},s-\frac{\lambda-t}{2}\right) ds\\
			&=\int_{-\infty}^{\infty} W(-s+d) f\left(s+\frac{\lambda-t}{2},s-\frac{\lambda-t}{2}\right) ds\qquad (\text{as } W(x)=W(-x),\; \forall\; x\in \Re)\\
			&=\int_{-\infty}^{\infty} W(s+d) f\left(-s+\frac{\lambda-t}{2},-s-\frac{\lambda-t}{2}\right) ds\\
			&=\int_{-\infty}^{\infty} W(s+d) f\left(s+\frac{\lambda-t}{2},s-\frac{\lambda-t}{2}\right) ds \qquad (\text{as } f(z_1,z_2)=f(-z_2,z_1),\; \forall \; (z_1,z_2)\in \Re^2)\\
			&=r^{*}_{\lambda}(-d,t).
		\end{align*}
		The assumption (D1) ensures that $r^{*}_{\lambda}(d,t)$ is a strictly convex function of $d$. Thus 
		\begin{align*}
			r^{*}_{\lambda}(0,t)&=r^{*}_{\lambda}\left(\frac{d+(-d)}{2},t\right)\\
			&\leq \frac{r^{*}_{\lambda}(d,t)+r^{*}_{\lambda}(-d,t)}{2}
			=r^{*}_{\lambda}(d,t)\\
			\implies r_{\lambda}\left(\frac{\lambda-t}{2},t\right)&\leq r_{\lambda}(c,t).
		\end{align*}
	\end{proof}

	The following theorem provide a sufficient condition under which a location equivariant estimator of $(\theta_1,\theta_2)$ is inadmissible. In such cases, the theorem also provides dominating estimators.
	
	\vspace*{2mm}

	\noindent \textbf{Theorem 2.1.1.} Under the assumptions (D1) and (D2), let $\boldsymbol{\delta}_{\psi}(\bold{X})=(X_1-\psi(D),X_2+\psi(D))$ be a location equivariant estimator of $(\theta_1,\theta_2)\in \Theta_0$, where $\psi:\Re\rightarrow \Re$. Suppose that $P_{\boldsymbol{\theta}}\left(\psi(D)<\frac{-D}{2}\right)>0,\;\forall \; \boldsymbol{\theta}\in \Theta_0$. Define $ \psi^{*}(t)=\max\{\frac{-t}{2},\psi(t)\},\; t\in\Re$. Then 
	$$R(\boldsymbol{\theta},\boldsymbol{\delta}_{\psi^{*}})\leq R(\boldsymbol{\theta},\boldsymbol{\delta}_{\psi}),\; \forall \; \boldsymbol{\theta}\in \Theta_0,$$
	where $\boldsymbol{\delta}_{\psi^{*}}(\bold{X})=(X_1-\psi^{*}(D),X_2+\psi^{*}(D))$.
	\begin{proof} Define $A=\{t:\psi(t)<\frac{-t}{2}\}$ and $B=\{t:\psi(t)\geq \frac{-t}{2}\}$, so that $\psi^{*}(t)=\frac{-t}{2}$, if $t\in A$, and $\psi^{*}(t)=\psi(t)$, if $t\in B$. The Lemma 2.1.1 and the assumption (D1) imply that, for any fixed $t\in \Re$ and $\lambda\geq 0$,
		$$r_{\lambda}(c,t)=\int_{-\infty}^{\infty} W(s-c)f(s,s+ t-\lambda)ds $$
		is decreasing in $c\in (-\infty,\frac{\lambda-t}{2}]$, increasing in $c\in [\frac{\lambda-t}{2},\infty),$ with unique minimum at $c\equiv \frac{\lambda-t}{2}$. Since, for any fixed $\lambda\geq 0$ and any $t$, $\frac{-t}{2}\leq \frac{\lambda-t}{2}< \infty$, it follows that, for any $\lambda\geq 0$, $r_{\lambda}(c,t)$
		is decreasing in $c\in (-\infty,\frac{-t}{2}]$.
		Consequently, for any $\lambda\geq 0$, $r_{\lambda}(\psi(t),t)\geq r_{\lambda}(\frac{-t}{2},t)$, for $t\in A$. Therefore,
		\begin{align*}
			R(\boldsymbol{\theta},\boldsymbol{\delta}_{\psi})&=2 \int_{-\infty}^{\infty} r_{\lambda}(\psi(t),t)\, dt\\
			&=2\left[\int_{A}\! r_{\lambda}(\psi(t),t)\,dt + \int_{B}\! r_{\lambda}(\psi(t),t)\, dt\right]\\
			&\geq 2\left[\int_{A} \!r_{\lambda}\left(\frac{-t}{2},t\right)\, dt + \int_{B}\! r_{\lambda}(\psi(t),t)\,dt\right]\\
			&=2\int_{-\infty}^{\infty} r_{\lambda}(\psi^{*}(t),t)\,0 dt\\
			&=R(\boldsymbol{\theta},\boldsymbol{\delta}_{\psi^{*}}), \; \; \boldsymbol{\theta}\in \Theta_0.
		\end{align*}
	\end{proof}

	\noindent 	The proof of the following Corollary is contained in the proof of Theorem 2.1.1, and hence omitted.
	\vspace*{2mm}
	
	\noindent	\textbf{Corollary 2.1.1.} Suppose that the assumptions (D1) and (D2) hold. Let $\boldsymbol{\delta}_{\psi}(\bold{X})=(X_1-\psi(D),X_2+\psi(D))$ be a location equivariant estimator of $\boldsymbol{\theta}$, where $\psi:\Re\rightarrow \Re$ is such that $P_{\boldsymbol{\theta}}\left(\psi(D)<\frac{-D}{2}\right)>0,\;\forall\; \boldsymbol{\theta}\in \Theta_0$. Let $\psi_0:\Re\rightarrow \Re$ be such that $\psi(t)\leq \psi_0(t) <\frac{-t}{2}$, whenever $\psi(t)<\frac{-t}{2}$, and $\psi_0(t)=\psi(t)$, whenever $\psi(t)\geq \frac{-t}{2} $. Then, 
	$R(\boldsymbol{\theta},\boldsymbol{\delta}_{\psi_0})$ $\leq R(\boldsymbol{\theta},\boldsymbol{\delta}_{\psi}),\; \forall \; \boldsymbol{\theta}\in \Theta_0,$
	where $\boldsymbol{\delta}_{\psi_0}(\bold{X})=(X_1-\psi_0(D),X_2+\psi_0(D))$.
	
	\vspace*{4mm}

	Under the unrestricted parameter space $\Theta=\Re^2$, it is easy to verify that the unrestricted best location equivariant estimator (BLEE) of $\boldsymbol{\theta}$ is
	$\boldsymbol{\delta}_0(\bold{X})$ $=(X_1,X_2)$. Using Theorem 2.1.1, we conclude that the unrestricted BLEE $\delta_0(\bold{X})$ is inadmissible  for estimating $\boldsymbol{\theta}$ and is dominated by the estimator $$\boldsymbol{\delta}_{\psi_0^{*}}(\bold{X})=(X_1-\psi_0^{*}(D),X_2+\psi_0^{*}(D))=\begin{cases}
		(X_1,X_2),&\text{if } X_1\leq X_2\\
		\left(\frac{X_1+X_2}{2},\frac{X_1+X_2}{2}\right),&\text{if }\; X_1>X_2
	\end{cases},$$
	where $\psi_0^{*}(D) =\max\big\{\frac{-D}{2},0\big\}$.
	
	\vspace*{3mm}

	Now we consider isotonic regression estimators (or mixed estimators) based on the unrestricted BLEE 
	$\boldsymbol{\delta}_0(\bold{X})$ $=(X_1,X_2)$. Let $\mathcal{D}=\{\boldsymbol{\delta}_{\alpha}:-\infty<\alpha<\infty\}$ be the class of isotonic regression estimators of $(\theta_1,\theta_2)$ based on the unrestricted BLEE $\boldsymbol{\delta}_0(\bold{X})=(X_1,X_2)$, where
	\begin{align*}
		\boldsymbol{\delta}_{\alpha}(\bold{X})&=(\delta_{1,\alpha}(\bold{X}),\delta_{2,\alpha}(\bold{X}))\\
		&= \begin{cases}
			(X_1,X_2), &\text{if } X_1\leq X_2\\
			(\alpha X_1 +(1-\alpha) X_2, (1-\alpha)X_1 +\alpha X_2), &\text{if } X_1>X_2
		\end{cases}\\
		&=(X_1-\psi_{\alpha}(D),X_2+\psi_{\alpha}(D)),
	\end{align*}
	where $\psi_{\alpha}(t)=\begin{cases} 0, &\text{if } t\geq 0 \\
		-(1-\alpha)t, &\text{if } t<0 \end{cases}.$
	
	\vspace*{2mm}
	
	Note that $P_{\boldsymbol{\theta}}\left(\psi_{\alpha}(D)< \frac{-D}{2}\right)>0,\;\forall\; \boldsymbol{\theta}\in \Theta_0$, if, and only if, $\alpha\geq \frac{1}{2}.$ Using Corollary 2.1.1, we conclude that, for $\frac{1}{2}\leq \alpha_0<\alpha_1<\infty$, the estimator $\boldsymbol{\delta}_{\alpha_0}$ dominates the estimator $\boldsymbol{\delta}_{\alpha_1}$. For the independent and bivariate normal probability model and the sum of squared errors loss function, the above consequences of Theorem 2.1.1 and Corollary 2.1.1 are obtained in Kumar and Sharma (1988) and Patra and Kumar (2017).
	\vspace*{3mm}


	\subsection{A Class of Improved Estimators Over the BLEE $(X_1,X_2)$}\label{sec2.2}
	\vspace*{2mm}
	\noindent
	
	\vspace*{2mm}
	
	In this section, we apply the Kubokawa (1994) technique to obtain a class of estimators improving over the BLEE $(X_1, X_2) $. Further, we obtain the Brewster-Zidek (1974) type and the Stein (1964) type improved estimators over the BLEE.
	
	Consider estimation of $(\theta_1,\theta_2)$ under the loss function \eqref{eq:2.2}, when it is known apriori that $\boldsymbol{\theta}\in\Theta_0$. Throughout this section, we will assume that the function $W(\cdot)$ is absolutely continuous and satisfies the assumption (D2). 
	\vspace*{2mm}

	The following lemma will be useful in proving the main results of this section. The proof of the lemma is straight forward and hence omitted. 
	\\~\\ \textbf{Lemma 2.2.1.} Let $s_0\in \Re$ and let $M:\Re\rightarrow\Re$ be such that $M(s)\leq 0,\; \forall \; s<s_0, $ and $M(s)\geq 0,\; \forall \; s> s_0$. Let $ M_i:\Re\rightarrow [0,\infty), \; i=1,2,$ be non-negative functions such that
	$M_1(s) M_2(s_0) \geq (\leq)\, M_1(s_0) M_2(s),\; \forall \; s<s_0,
	\text{ and } M_1(s) M_2(s_0) \leq\,(\geq)\; M_1(s_0) M_2(s),\; \forall \; s$ $>s_0.$
	Then, 
	$$ M_2(s_0) \int\limits_{-\infty}^{\infty} M(s) \, M_1(s) ds\leq\;(\geq)\; M_1(s_0) \int\limits_{-\infty}^{\infty} M(s) \, M_2(s) ds.$$
	
	The facts stated in the following lemma are well known in the theory of stochastic orders (see Shaked and Shanthikumar (\citeyear{MR2265633})). The proof of the lemma is straight forward, hence skipped.
	\\~\\ \textbf{Lemma 2.2.2.} If, for any fixed $\Delta\geq 0$ and $t\in \Re$, $\frac{f(s,s+t-\Delta)}{f(s,s+t)}$ is increasing (decreasing) in $s$,
	then $\frac{\int_{-\infty}^{t-\Delta}f(s,s+y)dy}{\int_{-\infty}^{t}f(s,s+y)dy}$ is increasing (decreasing) in $s$ and $\frac{f(s,s+t)}{\int_{-\infty}^{t}f(s,s+y)dy}$ is decreasing (increasing) in $s.$
	
	\vspace*{2mm}

	\noindent	In the following theorem, we provide a class of estimators that improve upon the BLEE $\boldsymbol{\delta}_{0}(\boldsymbol{X})=(X_1,X_2)$.
	\\~\\ \textbf{Theorem 2.2.1.} Suppose that the assumptions (D1) and (D2) hold. Let $\boldsymbol{\delta}_{\psi}(\boldsymbol{X})=(X_1 - \psi(D),X_2+\psi(D))$ be a location equivariant estimator of $(\theta_1,\theta_2)$ such that $\psi(t)$ is decreasing (increasing) in $t$, $\lim_{t\to\infty} \psi(t)=0$ and $\int_{-\infty}^{\infty} \int_{-\infty}^{t} W^{'}(s-\psi(t))\; f(s,s+y)\,dy\,ds \, \geq \,(\leq) \,0,\; \forall\; t.$ Then 
	$$R(\boldsymbol{\theta},\boldsymbol{\delta}_{\psi})\leq R(\boldsymbol{\theta},\boldsymbol{\delta}_{0}), \;\;\; \forall \; \; \boldsymbol{\theta} \in \Theta_0.$$
	\begin{proof}
		Let us fix $\boldsymbol{\theta}\in\Theta_0$ and let $\lambda=\theta_2-\theta_1$, so that $\lambda \geq 0$. Let $Z_i=X_i-\theta_i,\; i=1,2,$ and $Z=Z_2-Z_1$. Consider the risk difference
		\begin{align*}
			\Delta(\lambda)&= R(\boldsymbol{\theta},\boldsymbol{\delta}_{0})-R(\boldsymbol{\theta},\boldsymbol{\delta}_{\psi})\\
			& =2 E_{\boldsymbol{\theta}}[W(Z_1)- W(Z_1-\psi(Z+\lambda))] \\
			&= 2E_{\boldsymbol{\theta}}\left[\int_{Z+\lambda}^{\infty}\Big\{ \frac{d}{dt} W(Z_1-\psi(t))\Big\}\; dt\right]\\
			&= 2 E_{\boldsymbol{\theta}}\left[\int_{Z}^{\infty} (-\psi^{'}(t+\lambda)) W^{'}(Z_1-\psi(t+\lambda))\; dt\right]\\
			&= -2\int_{-\infty}^{\infty} \psi^{'}(t+\lambda) E_{\boldsymbol{\theta}}[W^{'}(Z_1-\psi(t+\lambda))\; I_{(-\infty,t]}(Z)\;]\,dt,
		\end{align*}
		where, for any set $A$, $I_A(\cdot)$ denotes its indicator function. Since $\psi(t)$ is a decreasing (increasing) function of $t$, it suffices to show that, for every $t$ and $\lambda\geq 0$, 
		\begin{align}\label{eq:2.1.1}	E_{\boldsymbol{\theta}}[W^{'}(Z_1-\psi(t+\lambda))\; I_{(-\infty,t]}(Z)\;]\geq \, (\leq)\,0.
		\end{align}

		Since $W^{'}(t)$ is an increasing function of $t$ and $\psi(t)$ is a decreasing (increasing) function of $t$, for $\lambda\geq 0$, we have 
		\begin{align*}
			E_{\boldsymbol{\theta}}[W^{'}(Z_1-\psi(t+\lambda))\; I_{(-\infty,t]}(Z)\;]
			&\geq \,(\leq)\,	E_{\boldsymbol{\theta}}[W^{'}(Z_1-\psi(t))\; I_{(-\infty,t]}(Z)\;]\nonumber\\
			&=\int_{-\infty}^{\infty}\int_{-\infty}^{t} W^{'}(s-\psi(t))\;f(s,s+y)\,dy\,ds
		\end{align*}
		which, in turn, implies \eqref{eq:2.1.1}.
	\end{proof}

	\noindent	Now we will prove a useful corollary to the above theorem. The following corollary provides the Brewster-Zidek (1974) type (B-Z type) improvement over the BLEE $\delta_{0}(\bold{X})=(X_1,X_2)$.
	\\~\\ \textbf{Corollary 2.2.1. (i)}  Suppose that, for any fixed $\Delta\geq 0$ and $t$, $\frac{\int_{-\infty}^{t-\Delta}f(s,s+y)dy}{\int_{-\infty}^{t}f(s,s+y)dy}$ is increasing (decreasing) in $s$. Further suppose that, for every fixed $t$, the equation
	$$k_1(c\vert  t)=\int_{-\infty}^{\infty}\int_{-\infty}^{t} \; W^{'}(s-c)\; f(s,s+y)\,dy\,ds =0$$
	has the unique solution $c\equiv \psi_{0,1}(t)$. Then
	$$R(\boldsymbol{\theta},\boldsymbol{\delta}_{\psi_{0,1}})\leq R(\boldsymbol{\theta},\boldsymbol{\delta}_{0}), \;\;\; \forall \; \; \boldsymbol{\theta} \in \Theta_0,$$
	where $\boldsymbol{\delta}_{\psi_{0,1}}(\boldsymbol{X})=(X_1-\psi_{0,1}(D),X_2+\psi_{0,1}(D))$.
	\\~\\\textbf{(ii)} In addition to assumptions of (i) above, suppose that $\psi_{1,1}:\Re \rightarrow \Re $ is such that $ \psi_{1,1}(t) \leq \, (\geq) \, \psi_{0,1}(t), \; \forall \; t, \; \psi_{1,1}(t) $ is decreasing (increasing) in $t$ and $\lim_{t \to \infty}\; \psi_{1,1}(t) = 0$.
	Then 
	$$R(\boldsymbol{\theta}, \boldsymbol{\delta}_{\psi_{1,1}}) \leq R(\boldsymbol{\theta}, \boldsymbol{\delta}_{0}),\; \; \forall \; \boldsymbol{\theta} \in \Theta_0,$$
	where $\boldsymbol{\delta}_{\psi_{1,1}}(\boldsymbol{X})=(X_1-\psi_{1,1}(D),X_2+\psi_{1,1}(D))$.
	\begin{proof}
		It suffices to show that $\psi_{0,1}(t)$ satisfies conditions of Theorem 2.2.1. Note that the hypothesis of the corollary ensure that $\lim_{t\to\infty} \psi_{0,1}(t)=0$. To show that $\psi_{0,1}(t)$ is a decreasing (increasing) function of $t$, suppose that, there exist numbers $t_1$ and $t_2$ such that $t_1<t_2$ and $\psi_{0,1}(t_1)\neq \psi_{0,1}(t_2).$ We have $k_1(\psi_{0,1}(t_1)\vert t_1)=0$. Also, using the hypotheses of the corollary and the assumption (D1), it follows that $\psi_{0,1}(t_2)$ is the unique solution of $k_1(c\vert t_2)=0$ and $k_1(c\vert t_2)$ is a decreasing function of $c$. Let $ s_0 = \psi_{0,1}(t_1), \; M(s)=W^{'}(s-s_0),\; M_1(s)=\int_{-\infty}^{t_2}f(s,s+y)dy$ and $M_2(s)=\int_{-\infty}^{t_1}f(s,s+y)dy$. Then, under assumption (D1), using Lemma 2.2.1, we get
		\begin{small}
			$$\int_{-\infty}^{t_1}f(\psi_{0,1}(t_1),\psi_{0,1}(t_1)+y)\,dy\; \left(\int_{-\infty}^{\infty}\; W^{'}(s-\psi_{0,1}(t_1))\;\int_{-\infty}^{t_2} f(s,s+y)\,dy \;ds\right) \qquad \qquad \qquad \qquad \qquad \qquad \quad$$   $$ \qquad \leq\; (\geq)\; \int_{-\infty}^{t_2}f(\psi_{0,1}(t_1),\psi_{0,1}(t_1)+y)dy\; \left( \int_{-\infty}^{\infty}\; W^{'}(s-\psi_{0,1}(t_1))\;\int_{-\infty}^{t_1} f(s,s+y)\,dy \;ds \right)=0.$$
		\end{small}
		$$\implies \;\; k_1(\psi_{0,1}(t_1)\vert t_2)=\int_{-\infty}^{\infty}\int_{-\infty}^{t_2}\; W^{'}(s-\psi_{0,1}(t_1))f(s,s+y)\,dy\, ds\leq \; (\geq) \; 0.$$
		This implies that $k_1(\psi_{0,1}(t_1)\vert t_2)\,< \, (>) \, 0$, as $k_1(c\vert t_2)=0$ has the unique solution $c\equiv \psi_{0,1}(t_2)$ and $\psi_{0,1}(t_1)\neq\psi_{0,1}(t_2)$. Since $k_1(c\vert t_2)$ is a decreasing function of c, $k_1(\psi_{0,1}(t_2)\vert t_2)$ $=0$ and $k_1(\psi_{0,1}(t_1)\vert t_2)\,< \, (>) \, 0$, it follows that $\psi_{0,1}(t_1)>(<) \psi_{0,1}(t_2)$. 
		\\~\\ The proof of part (ii) is an immediate by-product of Theorem 2.2.1 using the fact that, for any $t$, $k_1(c\vert t)$ is a decreasing function of $c\in\Re$.
	\end{proof}

	\noindent \textbf{Remark 2.2.1.}	It is straightforward to see that the Brewster-Zidek (1974) type estimator $\delta_{\psi_{0,1}}(\cdot)$, derived in Corollary 2.2.1 (i), is the generalized Bayes estimator with respect to the non-informative prior $\pi(\theta_1,\theta_2)=1,\;(\theta_1,\theta_2)\in\Theta_0.$\vspace*{2mm}

	In the following section, we will provide an application of results derived in Sections 2.1-2.2 and validate the results through a simulation study.
	\section{\textbf{An Application and a Simulation study: Bivariate Normal Distribution}}\label{sec3}
	
	\noindent
	
	\vspace*{3mm}

	\noindent   Let $\bold{X}=(X_1,X_2)\sim BVN(\theta_1,\theta_2,\sigma^2,\sigma^2,\rho)$, where $(\theta_1,\theta_2)\in \Theta_0$ is the vector of unknown means, $\sigma>0$ is the common known standard deviation and $\rho \in (-1,1)$ is the known correlation coefficient. The joint pdf of $(Z_1,Z_2)=(X_1-\theta_1,X_2-\theta_2)$ is
	$$ f(z_1,z_2) =\frac{1}{2 \pi \sigma^2 \sqrt{1-\rho^2}} e^{-\frac{1}{2(1-\rho^2)\sigma^2}\left[z_1^2-2 \rho \, z_1 z_2+z_2^2\right]},\; \; \; \bold{z}=(z_1,z_2)\in \Re^2.$$ 
	Consider estimation of $(\theta_1,\theta_2)$ under the general loss function  
	\begin{equation}\label{eq:3.1}
		L(\boldsymbol{\theta},\bold{a})= W(a_1-\theta_1) + W(a_2-\theta_2),\;\boldsymbol{\theta}=(\theta_1,\theta_2)\in\Theta_0,\;\bold{a}=(a_1,a_2)\in\Re^2,
	\end{equation}
	where $W:\Re \rightarrow \Re$ is such that the assumption (D1) holds.\vspace*{2mm}

	Let $\boldsymbol{\delta}_{\psi}(\bold{X})=(X_1-\psi(D),X_2+\psi(D))$ be a location equivariant estimator of $\boldsymbol{\theta}$ and let $\psi^{*}(t)=\max\{\psi(t),\frac{-t}{2}\},\;t\in \Re$ be as defined in Theorem 2.1.1. Using Theorem 2.1.1,
	it follows that, if $P_{\boldsymbol{\theta}}\left[\psi(D)< \frac{-D}{2}\right]>0,\; \forall \;\boldsymbol{\theta}\in\Theta_0$, then the estimator $\delta_{\psi}(\bold{X})$ is inadmissible for estimating $\boldsymbol{\theta}$ and is dominated by $\boldsymbol{\delta}_{\psi^{*}}(D)=(X_1-\psi^{*}(D),X_2+\psi^{*}(D))$. 
	\vspace*{2mm}
	
	The unrestricted BLEE of $\boldsymbol{\theta}$ is $\boldsymbol{\delta}_0(\bold{X})=(X_1,X_2)$. Then, the BLEE $(X_1,X_2)$ is improved on by 
	\begin{align}\label{eq:3.2}
		\boldsymbol{\delta}_{RMLE}(\bold{X})&=\left(X_1-\max\Big\{0, \frac{-D}{2}\Big\},X_2+\max\Big\{0, \frac{-D}{2}\Big\}\right) \nonumber\\
		&= \left(\min\Big\{X_1, \frac{X_1+X_2}{2}\Big\},\max\Big\{X_2, \frac{X_1+X_2}{2}\Big\}\right).
	\end{align}
	It is easy to verify that $\boldsymbol{\delta}_{RMLE}$ is the restricted maximum likelihood estimator (MLE) of $\boldsymbol{\theta}$ under the restricted parameter space $\Theta_0$ (see Kumar and Sharma (\citeyear{MR981031}) and Patra and Kumar (\citeyear{patra2017estimating})). \vspace*{2mm}
	
	When $W(t)=t^2,\;t\in \Re,$ using Corollary 2.2.1, under the loss function \eqref{eq:3.1}, the Brewster-Zidek (1974) type (B-Z type) improvements over the BLEE $(X_1,X_2)$ is
	\begin{align}\label{eq:3.3} 
		\boldsymbol{\delta}_{\psi_{0,1}}(\boldsymbol{X})=\left(X_1-\frac{\tau}{2}\, \frac{\phi\left(\frac{D}{\tau}\right)}{\Phi\left(\frac{D}{\tau}\right)},X_2+\frac{\tau}{2}\, \frac{\phi\left(\frac{D}{\tau}\right)}{\Phi\left(\frac{D}{\tau}\right)}\right),
	\end{align}
	where $\tau=\sigma\sqrt{2(1-\rho)}$ and, $\phi(\cdot)$ and $\Phi(\cdot)$ are the p.d.f. and the d.f. of the standard normal distribution, respectively. When $W(t)=\vert t\vert,\;t\in \Re,$ using Corollary 2.2.1, under the loss function \eqref{eq:3.1}, the Brewster-Zidek (1974) type (B-Z type) improvements over the BLEE $(X_1,X_2)$ is
	\begin{align}\label{eq:3.4} 
		\boldsymbol{\delta}_{\psi_{0,1}}(\boldsymbol{X})=\left(X_1-C(D),X_2+C(D)\right),
	\end{align}
	where $C\equiv C(D)$ is the solution of the following equation
	\begin{equation}
		\int_{-\infty}^{C} \Phi\left(\frac{D+s(1-\rho)}{\sigma \sqrt{1-\rho^2}}\right) \phi\left(\frac{s}{\sigma}\right)ds=\frac{\sigma}{2}\Phi\left(\frac{D}{\tau}\right).
	\end{equation}
	Note that the estimators, given by \eqref{eq:3.3} and \eqref{eq:3.4}, are the generalized Bayes estimators with respect to the non-informative prior density (the Lebesgue measure) on $\Theta_0$ and the loss function \eqref{eq:3.1}, with $W(t)=t^2,\;t\in \Re,$ and $W(t)=\vert t\vert,\;t\in \Re$, respectively.

	\vspace*{6mm}

	\noindent \underline{\textbf{Simulation Study:}}

	\vspace*{2mm}

	For the above bivariate normal model, we considered estimation of vector $\boldsymbol{\theta}=(\theta_1,\theta_2)$ of means when it is known apriori that they satisfy the order restriction $\theta_1\leq \theta_2$. For estimation of $\boldsymbol{\theta}$, we obtained estimators (given by \eqref{eq:3.2}, \eqref{eq:3.3} and \eqref{eq:3.4}) improving on the BLEE $(X_1,X_2)$. The improved estimator \eqref{eq:3.2} is the same as the restricted maximum likelihood estimator (MLE), and the improved estimators \eqref{eq:3.3} and \eqref{eq:3.4} are the generalized Bayes (GB) estimators with respect to the squared error loss and the absolute error loss, respectively. To further evaluate the performances of the improved estimators, in this section, we compare the risk performances of the BLEE $(X_1,X_2)$, the restricted MLE (as defined in \eqref{eq:3.2}) and the generalized Bayes (GB) estimators (as defined in \eqref{eq:3.3} and \eqref{eq:3.4}), numerically, through the Monte Carlo simulations. For simulation study, we take $W(t)=t^2,\; t\in \Re$ (i.e., sum of squared error losses) and $W(t)=\vert t\vert ,\; t\in \Re$ (i.e., sum of absolute error losses). The simulated risks of the BLEE, the restricted MLE and the GB estimator have been computed.
	\vspace*{2mm}

\begin{figure}
		\begin{subfigure}{.48\textwidth}
			\centering
			\includegraphics[width=72mm,scale=1.2]{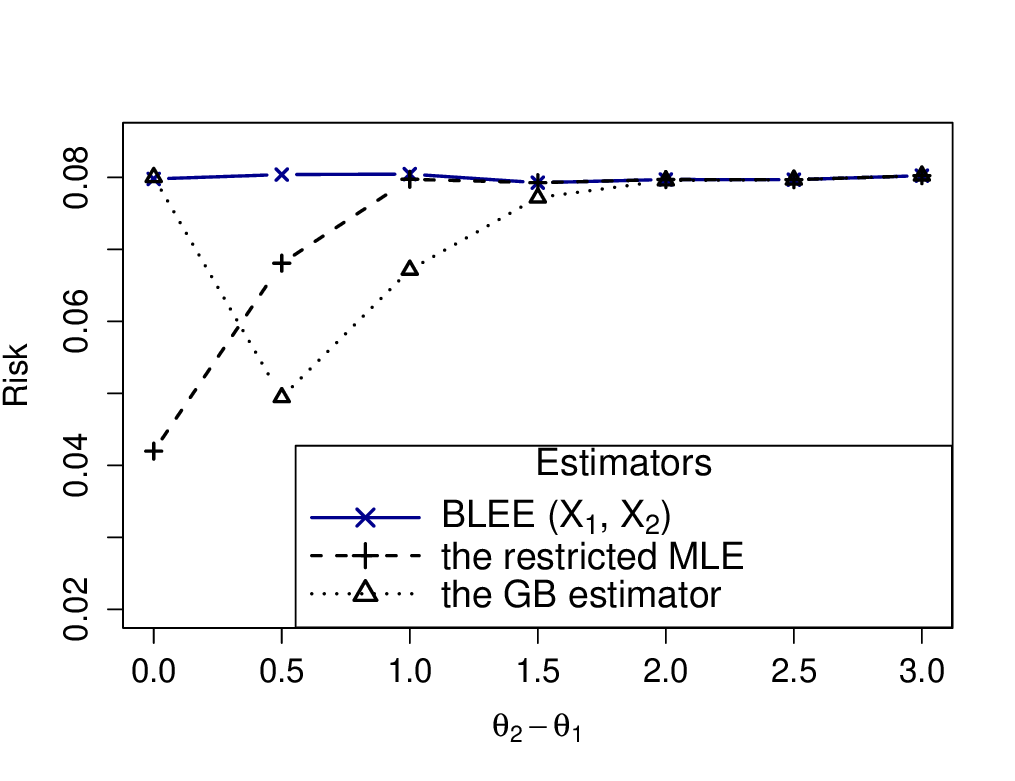} 
			\caption{$\sigma=0.2$ and $\rho=-0.9$.} 
			
		\end{subfigure}
		\begin{subfigure}{.48\textwidth}
			\centering
			\includegraphics[width=72mm,scale=1.2]{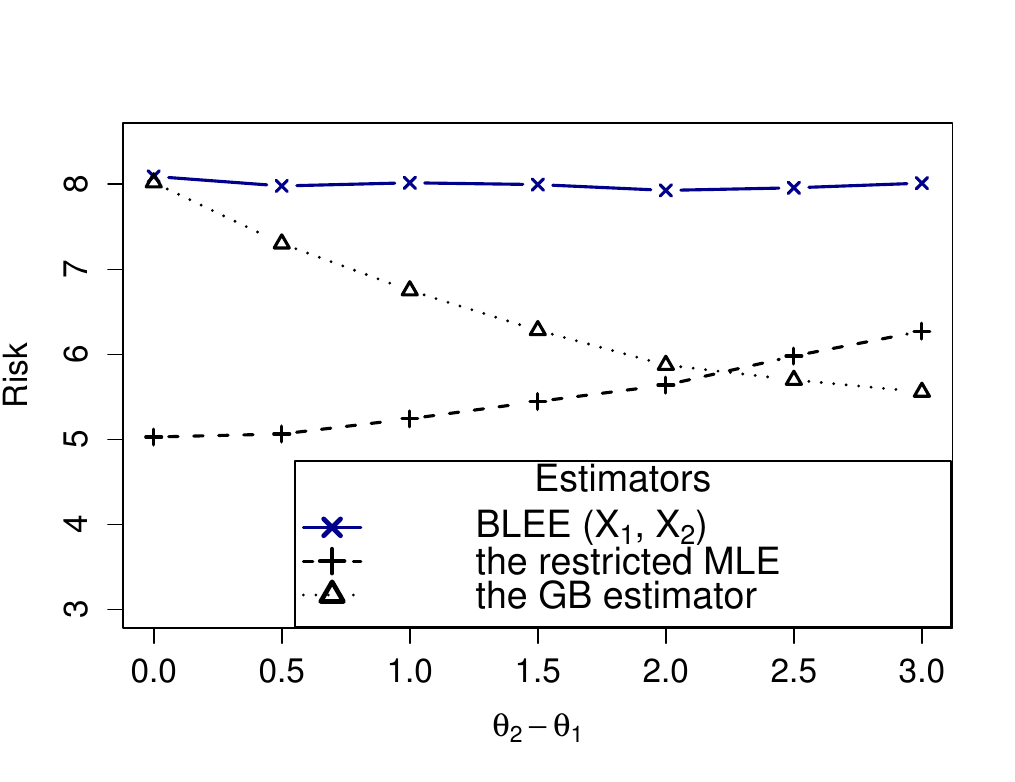} 
			
			\caption{$\sigma=2$ and $\rho=-0.5$.} 
		\end{subfigure}	\\	\begin{subfigure}{.48\textwidth}
			\centering
			\includegraphics[width=72mm,scale=1.2]{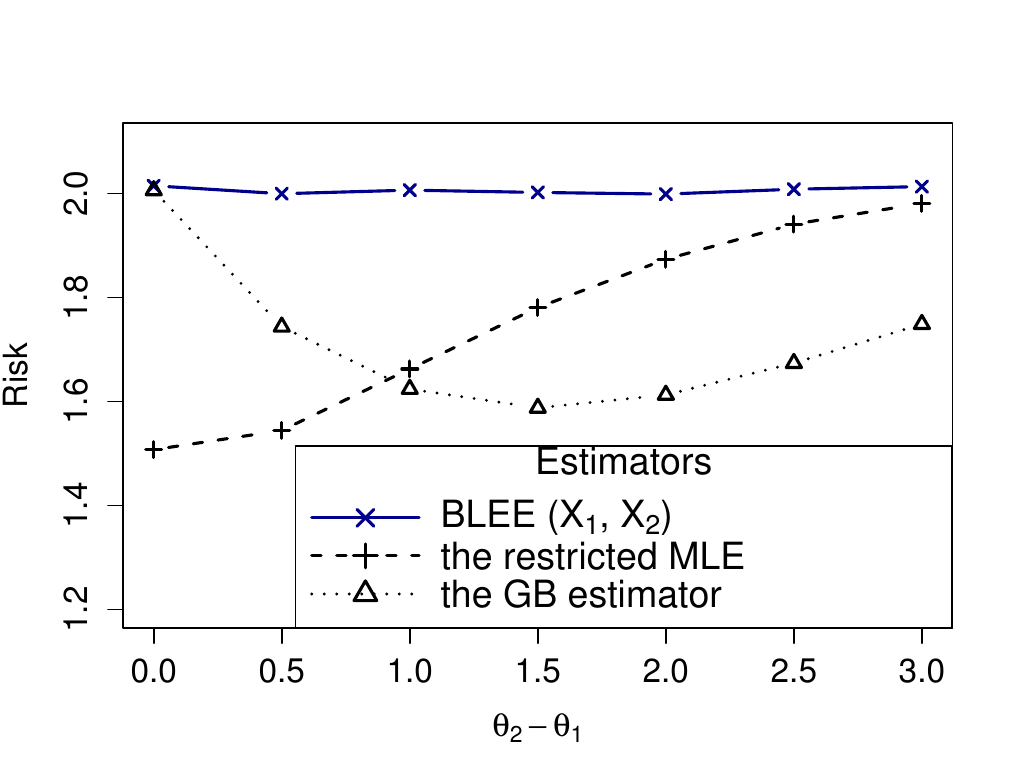} 
			\caption{$\sigma=1$ and $\rho=0$.} 
			
		\end{subfigure}
		\begin{subfigure}{.48\textwidth}
			\centering
			
			\includegraphics[width=72mm,scale=1.2]{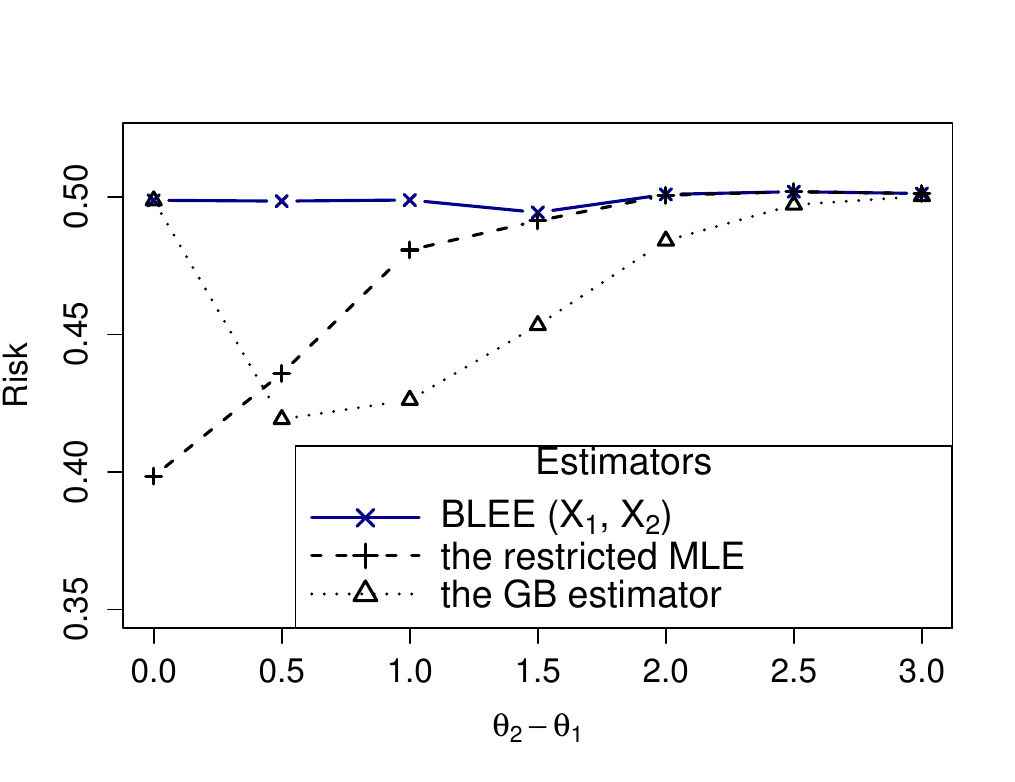} 
			
			\caption{$\sigma=0.5$ and $\rho=0.2$.} 
			
		\end{subfigure}
		\\	\begin{subfigure}{.48\textwidth}
			\centering
			\includegraphics[width=72mm,scale=1.2]{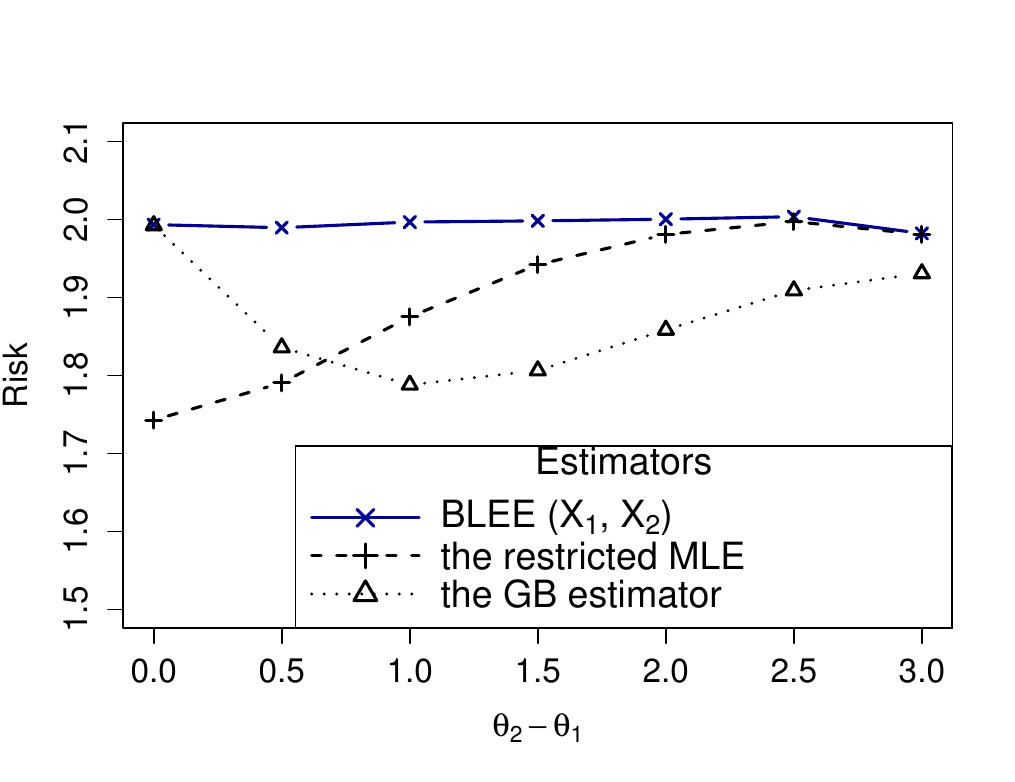} 
			\caption{ $\sigma=1$ and $\rho=0.5$.} 
			
		\end{subfigure}
		\begin{subfigure}{.48\textwidth}
			\centering
			
			\includegraphics[width=72mm,scale=1.2]{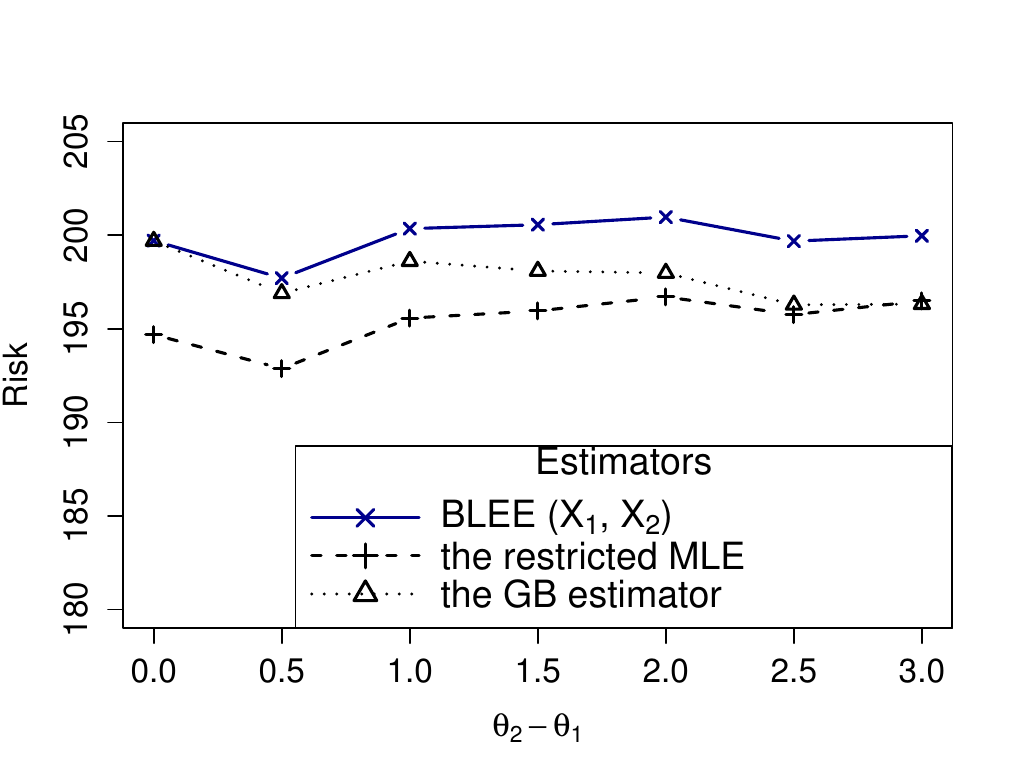} 
			\caption{$\sigma=10$ and $\rho=0.9$.} 
			
		\end{subfigure}
		
		\caption{Risk plots of estimators of location parameter $\boldsymbol{\theta}$ against values of $\theta_2-\theta_1$: when $W(t)=t^2,\;t\in \Re$.}
		\label{fig1}
	\end{figure}

\begin{figure}
		\begin{subfigure}{.48\textwidth}
			\centering
			\includegraphics[width=72mm,scale=1.2]{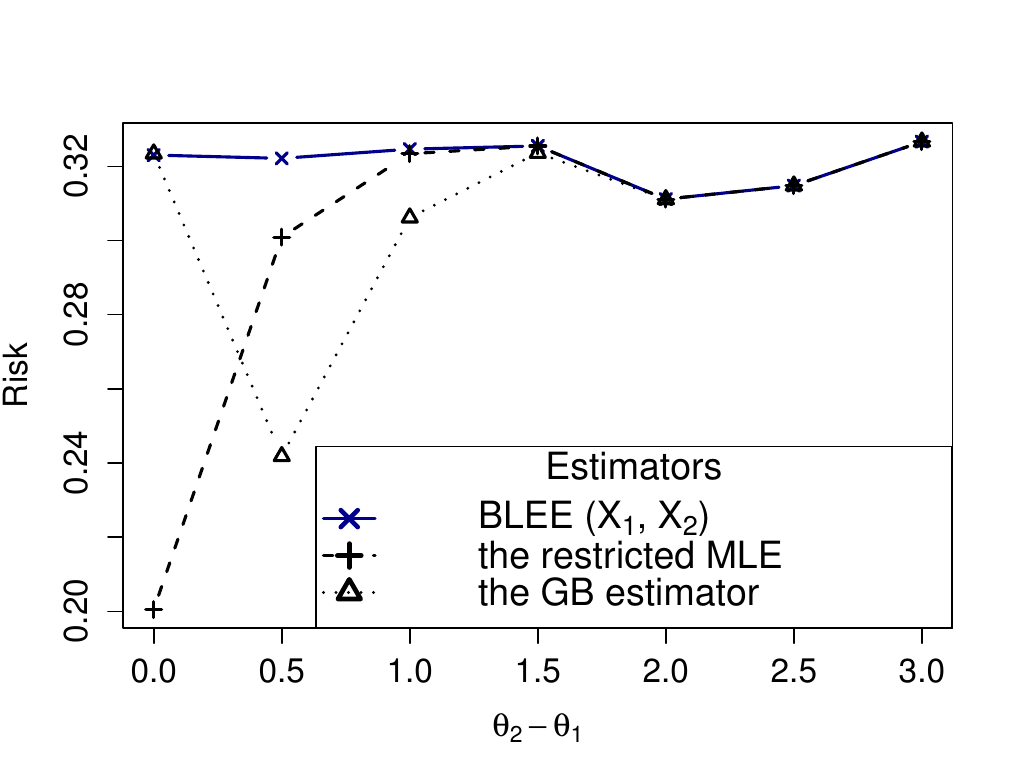} 
			\caption{$\sigma=0.2$ and $\rho=-0.9$.} 
			
		\end{subfigure}
		\begin{subfigure}{.48\textwidth}
			\centering
			\includegraphics[width=72mm,scale=1.2]{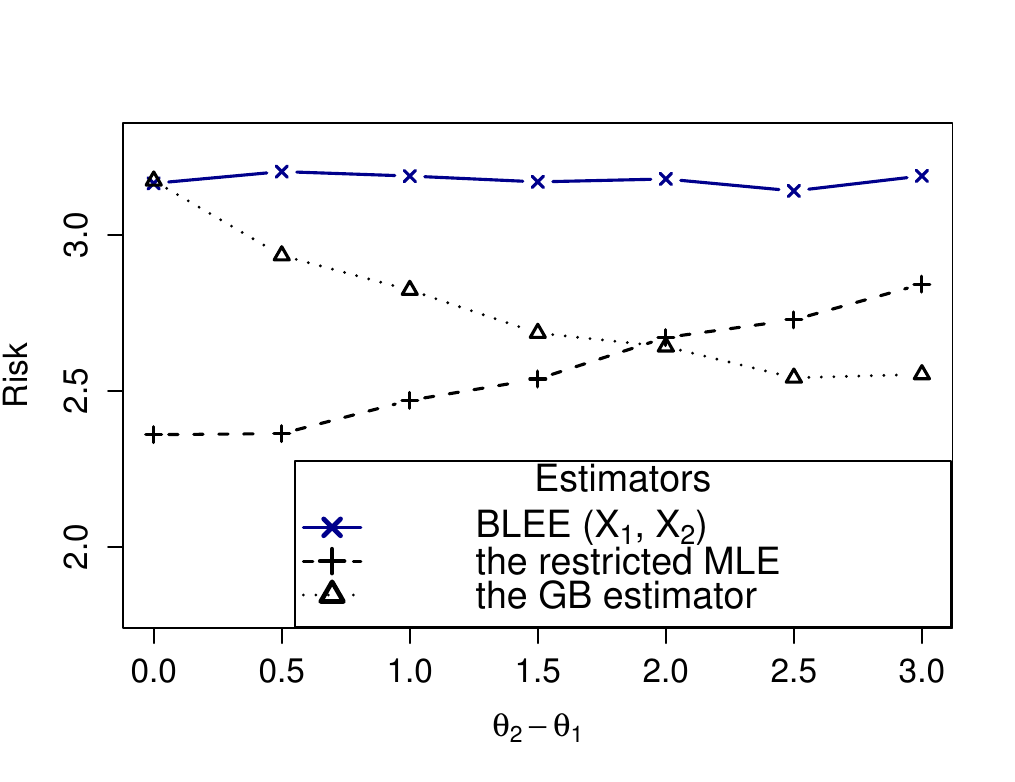} 
			
			\caption{$\sigma=2$ and $\rho=-0.5$.} 
			
		\end{subfigure}
		\\	\begin{subfigure}{.48\textwidth}
			\centering
			\includegraphics[width=72mm,scale=1.2]{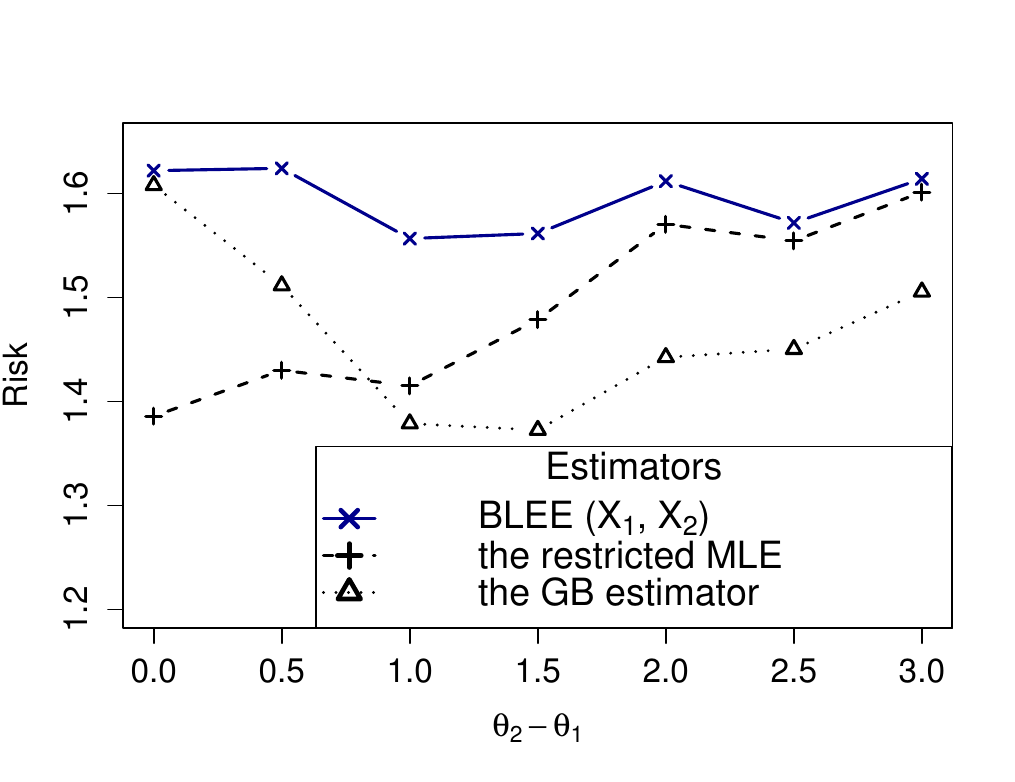} 
			\caption{$\sigma=1$ and $\rho=0$.} 
			
		\end{subfigure}
		\begin{subfigure}{.48\textwidth}
			\centering
			
			\includegraphics[width=72mm,scale=1.2]{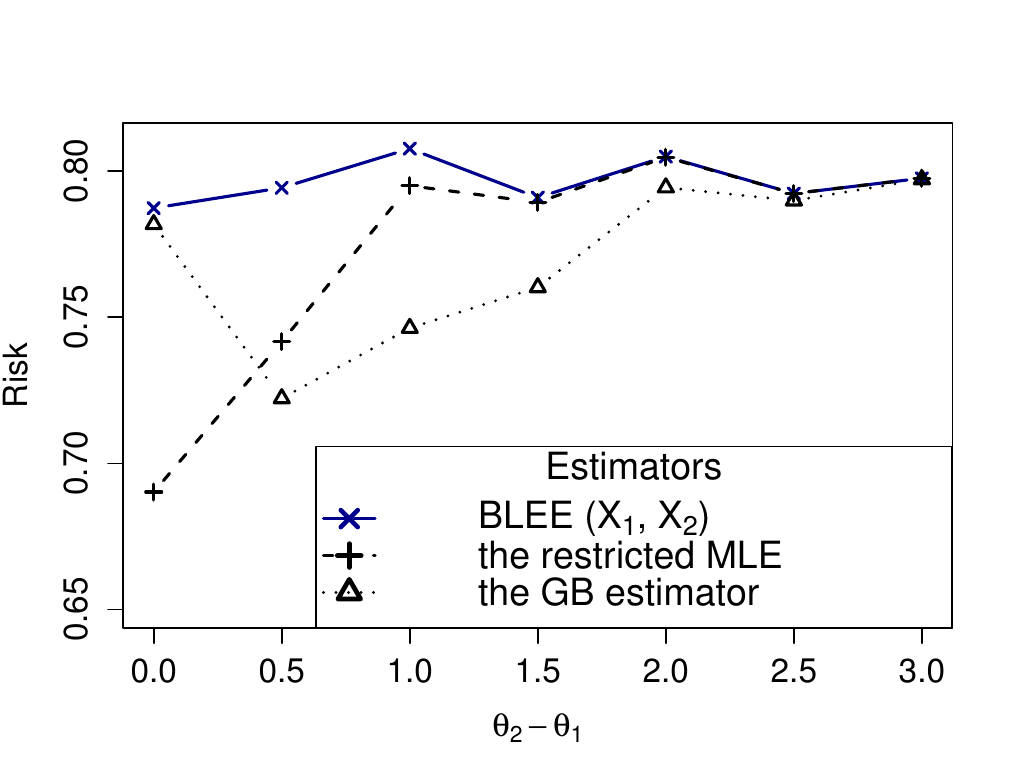} 
			
			\caption{$\sigma=0.5$ and $\rho=0.2$.} 
			
		\end{subfigure}
		\\	\begin{subfigure}{.48\textwidth}
			\centering
			\includegraphics[width=72mm,scale=1.2]{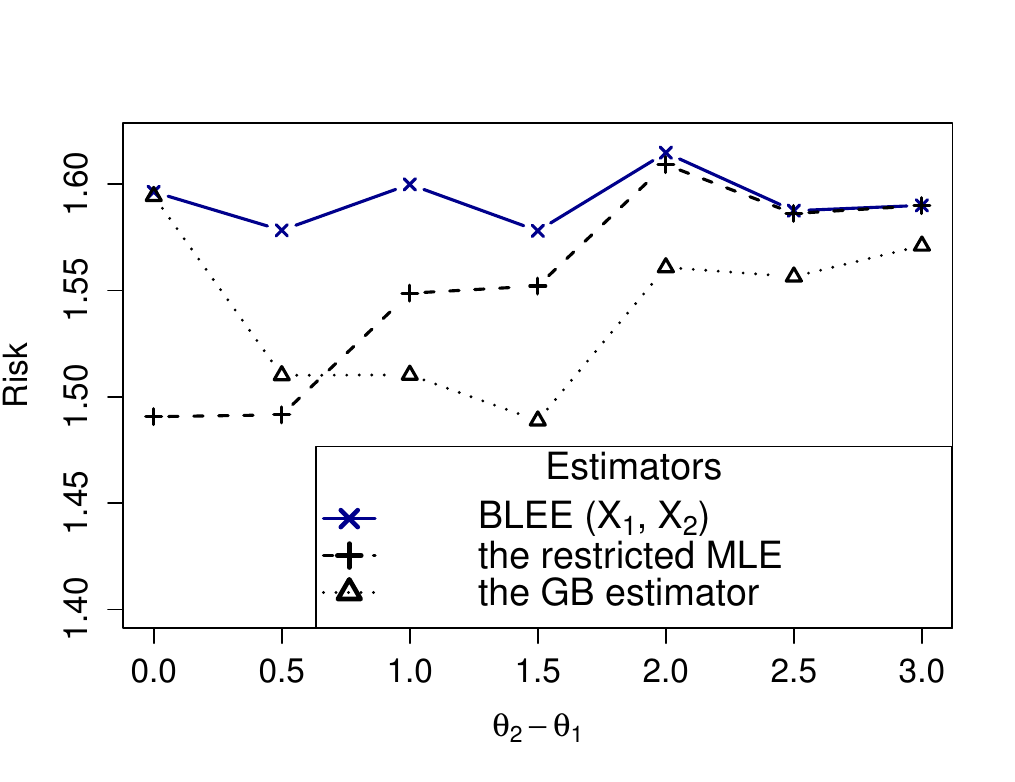} 
			\caption{ $\sigma=1$ and $\rho=0.5$.} 
			
		\end{subfigure}
		\begin{subfigure}{.48\textwidth}
			\centering
			
			\includegraphics[width=72mm,scale=1.2]{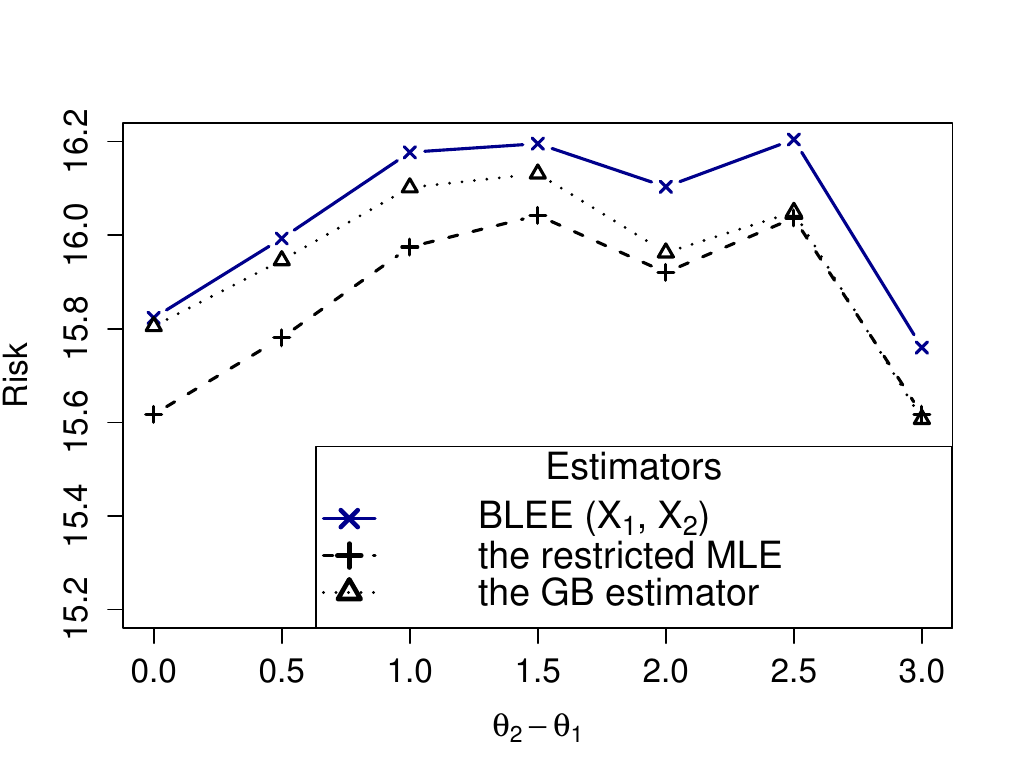} 
			\caption{$\sigma=10$ and $\rho=0.9$.} 
			
		\end{subfigure}
		
		\caption{Risk plots of estimators of location parameter $\boldsymbol{\theta}$ against values of $\theta_2-\theta_1$: when $W(t)=\vert t \vert,\;t\in \Re$.}
		\label{fig2}
	\end{figure}

	For simulations, 10,000 random samples of size 1 were generated from the relevant bivariate normal distribution. The simulated values of the risks of the BLEE, the restricted MLE and the GB estimator under the sum of squared error loss functions and the sum of absolute error loss functions are plotted in Figure \ref{fig1} and Figure \ref{fig2}, respectively. The following observations are evident from Figure \ref{fig1} and Figure \ref{fig2}:
	\vspace*{2mm}
	
	\noindent (i) The restricted MLE and the GB estimator perform better than the BLEE, which is in conformity with our theoretical findings;\vspace*{2mm}
	
	\noindent (ii) The performance of the restricted MLE is significantly better when $\theta_1$ and $\theta_2$ ($\theta_1\leq \theta_2$) are close, otherwise the GB estimator performs better.\vspace*{2mm}
	
	\noindent(iii) There is no clear cut winner between the restricted MLE and the GB estimator.\vspace*{2mm}
	
	\noindent(iii) Also, note that the performance of both the GB estimators (given by \eqref{eq:3.3} and \eqref{eq:3.4}) remain the same, relative to other two estimators, under both the loss functions, the squared error loss function and the absolute error loss function.

	\FloatBarrier
	\FloatBarrier

	\FloatBarrier
	\FloatBarrier

	\section{Real Life Data Analysis}\label{sec4}
	

	\begin{wraptable}{r}{4.5cm}
		\qquad	\begin{tabular}{| c|  c | c |} 
			\hline
			Girl/Boy & 8 year & 10 year \\ 
			[0.5ex] 
			\hline\hline
			Girl &	21 &	20 \\
			Girl &	21 &	21.5 \\
			Girl &	20.5 &	24 \\
			Girl &	23.5 &	24.5 \\
			Girl &	21.5 &	23 \\
			Girl &	20 &	21 \\
			Girl &	21.5 &	22.5 \\
			Girl &	23 &	23 \\
			Girl &	20 &	21 \\
			Girl &	16.5 &	19 \\
			Girl &	24.5 &	25 \\
			Boy &	26 &	25 \\
			Boy &	21.5 &	22.5 \\
			Boy	&	23 &	22.5 \\
			Boy &	25.5 &	27.5 \\
			Boy &	20 &	23.5 \\
			Boy &	24.5 &	25.5 \\
			Boy &	22 &	22 \\
			Boy &	24 &	21.5 \\
			Boy &	23 &	20.5 \\
			Boy &	27.5 &	28 \\
			Boy &	23 &	23 \\
			Boy &	21.5 &	23.5 \\
			Boy &	17 &	24.5 \\
			Boy &	22.5 &	25.5 \\
			Boy &	23 &	24.5 \\
			Boy &	22 &	21.5 \\
			\hline\hline
			Mean  & 22.185 & 23.167\\
			\hline
			
		\end{tabular}
		\caption{The size of pituitary fissure of children at different ages.}
		\label{table:1}
		
	\end{wraptable}

	We consider ``the dental study data", discussed by Potthoff and Roy (\citeyear{potthoff1964generalized}), that is presented in Table \ref{table:1}. This study was conducted at the University of North Carolina Dental School. In this study, the size (in millimeters) of the pituitary fissure was measured in children of different ages. To test the bivariate normality of this data set, we applied the Henze-Zirkler and Anderson-Darling multivariate normality tests and observed p-values of $0.474$ and $0.288$, respectively, suggesting that there is no significant departure from normality. Here, it is reasonable to assume that the size of the pituitary fissure increases with age. We performed the paired t-test and got the p-value $= 0.008$ in favour of the assumption that the pituitary fissure increases with age. We also performed the variance equality test of the data with respect to 8 year and 10 year and got the p-value $= 0.542$. As a result, we can say that the variances of both datasets are the same.
	\vspace*{3mm}

	Now, to illustrate the findings of our paper, suppose that the data of 5 girls and 8 boys, presented in Table \ref{table:2} is reported (for reference see p.2 of Robertson et al. (\citeyear{MR961262})). Let $X_1$ and $X_2$ be random variables representing the average size of the pituitary fissure of 8 year and 10 year children, respectively. Then $(X_1,X_2)$ follows a bivariate normal distribution with means $\theta_1$ and $\theta_2$, common variance $\sigma^2$ and correlation coefficient $\rho$. We use the common sample variance and the sample correlation of the data of Table \ref{table:1} as the plug-in values for $\sigma^2=\frac{5.435}{13}=0.418$ and $\rho=0.626$. We know that $\theta_1\leq \theta_2$. By exploiting this information, we obtain improvements over unrestricted estimators (based on the data in Table \ref{table:2}) for $\theta_1$ and $\theta_2$.

	\vspace*{2mm}
	In the starting of this section, we have seen that, under the general loss function \eqref{eq:3.1}, the improved BLEE (the Stein-type estimator \eqref{eq:3.2}) dominates the unrestricted BLEE $(X_1,X_2)$. The improved BLEE (the restricted MLE) of $(\theta_1,\theta_2)$ is 
	\begin{align*}\label{eq:3.2}
		\left(\min\Big\{X_1, \frac{X_1+X_2}{2}\Big\},\max\Big\{X_2, \frac{X_1+X_2}{2}\Big\}\right)=(22.86,22.86).
	\end{align*}
	Based on the theoretical results and the simulation studies, we infer that, rather than $23.077$ and $22.654$, respectively, $22.86$ be taken as a common estimated value for $\theta_1$ and $\theta_2$. Also, the Brewster-Zidek (1974) type improved estimated values (as defined by \eqref{eq:3.3} and \eqref{eq:3.4}) under the squared error loss and the absolute error loss are $(22.77,22.96)$ and $(22.71,23.03)$, respectively.

	\FloatBarrier
	\begin{table}[h!]
		\centering
		\begin{tabular}{| c|  c | c |} 
			\hline
			Girl/Boy & 8 year & 10 year \\ 
			[0.5ex] 
			\hline\hline
			Girl &	21 &	20 \\
			Girl &	21 &	21.5 \\
			Girl &	20 &	21 \\
			Girl &	23 &	23 \\
			Girl &	24.5 &	25 \\
			Boy &	26 &	25 \\
			Boy	&	23 &	22.5 \\
			Boy &	22 &	22 \\
			Boy &	24 &	21.5 \\
			Boy &	23 &	20.5 \\
			Boy &	27.5 &	28 \\
			Boy &	23 &	23 \\
			Boy &	22 &	21.5 \\
			\hline\hline
			Mean  & 23.077 & 22.654\\
			\hline
			
		\end{tabular}
		\caption{The size of pituitary fissure of children at different ages.}
		\label{table:2}
	\end{table}
	\FloatBarrier


	\section{Conclusions}\label{sec5}
	
	In this paper, we considered simultaneous estimation of order restricted location parameters of a general bivariate symmetric model, under a general loss function. We used the Stein technique to obtain truncated estimators and the Kubokawa (or IERD) technique to obtain a class of smooth estimators, which dominate the BLEE estimator. Additionally, we obtained the Brewster-Zidek smooth estimator, which is also the generalized Bayes estimator. Our findings demonstrate that the Stein-type estimator is robust, as it does not depend on the choice of the probability model or the loss function. We also conducted a simulation study to confirm the findings of the paper, and also provided a real-life application of the results obtained in the paper.
	\vspace*{2mm}
	
	One can think about extending the results of the paper from the general bivariate distribution to a general multivariate distribution. This seems to be a challenging problem, and it may be taken up in our future research.

	\subsection*{Financial disclosure}
	
	This work was supported by the Council of Scientific and Industrial Research (CSIR) under Grant [number 09/092(0986)/2018].
	
	\subsection*{Conflict of interest}
	
	There is no conflict of interest by authors.

	\bibliographystyle{apalike}
	
	\bibliography{references}

\end{document}